# Further Improvements on Waring's Problem

Li An-Ping

Beijing 100085, P.R. China
apli0001@sina.com

Abstract

In the paper, we will continue the investigation of Waring's problem, and give further improvements on the estimations of $G(k)$.*

* The article is first posted in the web of NSTL in 2019, https://preprint.nstl.gov.cn/preprint/

## 1. Introduction

To improve the estimate $G(k)$ of Waring's problem, in paper [1], we have introduced the method of parameterized recursions, and obtain substantive improvements over the earlier results of [3,5]. In this paper, with a modification of recursive process, we will give further improvements.

**Theorem 1.** For sufficient large $k$ ,

$$G(k) \le \begin{cases} 3.094686k + o(k), & \text{if } k \text{ is not a power of 2,} \\ 4k, & \text{otherwise.} \end{cases} \tag{1.1}$$

**Theorem 2.** For $5 \le k \le 20,$ let $F(k)$ be as in the List 1.1, then

$$G(k) \le F(k). \tag{1.2}$$

| $k$ | $F(k)$ | $k$ | $F(k)$ | $k$ | $F(k)$ | $k$ | $F(k)$ |
|---|---|---|---|---|---|---|---|
| 5 | 17 | 9 | 29 | 13 | 41 | 17 | 53 |
| 6 | 19 | 10 | 33 | 14 | 45 | 18 | 57 |
| 7 | 23 | 11 | 35 | 15 | 47 | 19 | 61 |
| 8 | 32 | 12 | 39 | 16 | 64 | 20 | 63 |

List 1.1

## 2. The Proofs of Theorems 1, 2.

Similarly, the Lemma 4.1 of [1] is still the outset of the paper, for the convenience, we repeat it in the following, the symbols and notations are same as before.

**Lemma 1.**

$$\begin{aligned} &J_0 \le ZPS_{s-1}(P) + J_1, \\ &J_i \ll U_i + V_i, \\ &U_i = S_{s-1}(P_i)^{1/2} Z_{i+1}^{(2s-3)/2} \left( P(\tilde{H}_i \tilde{Z}_i)^2 Z_{i+1} S_{s-1}(P_{i+1}) \right)^{1/2} \\ &V_i = S_{s-1}(P_i)^{1/2} Z_{i+1}^{(2s-3)/2} (\tilde{H}_i \tilde{Z}_i J_{i+1})^{1/2}, \qquad 1 \le i < k, \end{aligned} \tag{2.1}$$

Where $\tilde{H}_i = \prod_{j \le i} H_j, \tilde{Z}_i = \prod_{j \le i} Z_j$

Besides, we will require a known result in Waring's problem (see [2],[3])

**Lemma 2.** For $k > 2$, let $\varsigma(k) = 4k$ if $k$ is a power of $2$, or $3k/2$ otherwise. suppose that $s = 2t + v$, $t, v$ are two positive integers, if satisfying $s \geq \varsigma(k)$, and $v \cdot 2^{1-k} > \Delta(t)$, then $G(k) \leq s$.

Proof of Theorem 1:

Let

$$U_{i-1} = V_{i-1} \cdot H_i^{-\alpha/2} P^{-\tau_i/2}, \quad 1 \leq i \leq k. \tag{2.2}$$

The parameters $\alpha, \tau_i, 1 \leq i \leq k$, will be decided later. Hence, it has

$$\left( P(\tilde{H}_{i-1}\tilde{Z}_{i-1})^2 Z_i S_{s-1}(P_i) \right) = (\tilde{H}_{i-1}\tilde{Z}_{i-1} J_i) \cdot H_i^{-\alpha} P^{-\tau_i}.$$

And

$$J_i = P(\tilde{H}_{i-1}\tilde{Z}_{i-1}) Z_i S_{s-1}(P_i) H_i^{\alpha} P^{\tau_i}.$$

On the other hand,

$$J_i \approx U_i \cdot H_{i+1}^{\alpha/2} P^{\tau_{i+1}/2} = S_{s-1}(P_i)^{1/2} Z_{i+1}^{(2s-3)/2} \left( P(\tilde{H}_i \tilde{Z}_i)^2 Z_{i+1} S_{s-1}(P_{i+1}) \right)^{1/2} H_{i+1}^{\alpha/2} P^{\tau_{i+1}/2}$$

Combine the two equalities, it has

$$P(\tilde{H}_{i-1}\tilde{Z}_{i-1}) Z_i S_{s-1}(P_i) H_i^{\alpha} P^{\tau_i} = S_{s-1}(P_i)^{1/2} Z_{i+1}^{(2s-3)/2} \left( P(\tilde{H}_i \tilde{Z}_i)^2 Z_{i+1} S_{s-1}(P_{i+1}) \right)^{1/2} H_{i+1}^{\alpha/2} P^{\tau_{i+1}/2}$$

It follows,

$$\frac{S_{s-1}(P_i)^{1/2}}{S_{s-1}(P_{i+1})^{1/2}} = \frac{Z_{i+1}^{(2s-2)/2}(H_i) H_{i+1}^{\alpha/2} P^{\tau_{i+1}/2}}{P^{1/2} H_i^{\alpha} P^{\tau_i}},$$

and

$$\frac{S_{s-1}(P_i)}{S_{s-1}(P_{i+1})} = \frac{Z_{i+1}^{(2s-2)} (H_i)^{2-2\alpha} H_{i+1}^{\alpha}}{P^{1+2\tau_i - \tau_{i+1}}}.$$

Let $S_{s-1}(X) = X^{\lambda_{s-1}}$, it has

$$\frac{P_i^{\lambda_{s-1}}}{P_{i+1}^{\lambda_{s-1}}} = Z_{i+1}^{(2s-2)} \frac{(P/Z_i^k)^{2-2\alpha} (P/Z_{i+1}^k)^{\alpha}}{P^{1+2\tau_i - \tau_{i+1}}}.$$

i.e.

$$(Z_{i+1})^{\lambda_{s-1}} = Z_{i+1}^{(2s-2)} P^{1-\alpha-(2\tau_i - \tau_{i+1})} Z_i^{-k(2-2\alpha)} Z_{i+1}^{-k\alpha}.$$

and

$$\theta_i = \theta_{i+1} \frac{((2s-2) - \lambda_{s-1} - k\alpha)}{k(2-2\alpha)} + \frac{1}{2k} \frac{1-\alpha-(2\tau_i - \tau_{i+1})}{(1-\alpha)}. \tag{2.3}$$

Let

$$\eta_i = 2\tau_i - \tau_{i+1}, \quad \eta_i = \tau + \mu\sigma^{k-i}, \quad 1 \le i < k. \tag{2.4}$$

Then it is easy to deduce that

$$\tau_{i+1} = \tau - \mu\sigma^{k-i}\frac{1-(2\sigma)^i}{1-2\sigma}, \qquad (\tau = \tau_1). \tag{2.5}$$

Especially,

$$\tau_k = \tau - \mu\sigma\frac{1-(2\sigma)^{k-1}}{1-2\sigma}. \tag{2.6}$$

Denote by $a = \dfrac{(2s-2)-\alpha k-\lambda_{s-1}}{2k(1-\alpha)}$, $b_i = \dfrac{1}{2k}\dfrac{1-\alpha-\eta_i}{1-\alpha}$, by (2.3), it has

$$\theta_{k-i} = a^i\theta_k + \frac{1}{2k}\frac{1-a^i}{1-a} - \frac{1}{2k(1-\alpha)}\left(\tau\frac{1-a^i}{1-a} + \mu\sigma\frac{\sigma^i - a^i}{\sigma - a}\right). \tag{2.7}$$

Furthermore, by (2.2) with $i = k$, it has

$$\left(P(\tilde{H}_{k-1}\tilde{Z}_{k-1})^2 Z_k S_{s-1}(P_k)\right)^{1/2} = (\tilde{H}_{k-1}\tilde{Z}_{k-1}J_k)^{1/2} \cdot H_k^{-\alpha/2}P^{-\tau_k/2}.$$

And it is easy to know that $J_k = P\tilde{H}_k\tilde{Z}_k S_{s-1}(P_k)$, hence

$$P(\tilde{H}_{k-1}\tilde{Z}_{k-1})^2 Z_k S_{s-1}(P_k) = (\tilde{H}_{k-1}\tilde{Z}_{k-1}P\tilde{H}_k\tilde{Z}_k S_{s-1}(P_k)) \cdot H_k^{-\alpha}P^{-\tau_k}$$

$$\begin{aligned} &\Rightarrow 1 = H_k^{1-\alpha}P^{-\tau_k} \Rightarrow P^{\tau_k} = H_k^{1-\alpha} = (P/Z_k^k)^{1-\alpha} \Rightarrow Z_k^{k(1-\alpha)} = P^{1-\alpha-\tau_k} \\ &\Rightarrow \theta_k = \frac{1}{k}\cdot\frac{1-\alpha-\tau_k}{(1-\alpha)}. \end{aligned} \tag{2.8}$$

So, by (2.7) and (2.8), there is

$$\theta = \theta_1 = \frac{1-\alpha-\tau}{2k(1-\alpha)(1-a)} + \frac{1-\alpha-\tau}{2k(1-\alpha)}\cdot\left(\frac{1-2a}{(1-a)} + \bar{\mu}\zeta\right)a^{k-1}. \tag{2.9}$$

Where $\bar{\mu} = \mu/(1-\alpha-\tau)$, $\quad \zeta = \left(2\sigma\dfrac{1-(2\sigma)^{k-1}}{(1-2\sigma)} - \dfrac{\sigma}{a}\dfrac{(\sigma/a)^{k-1}-1}{(\sigma/a)-1}\right)$.

On the other hand, by Lemma 2.1 of [1] with $i = 1$, it has

$$S_s(\tilde{P}) \ll Z^{2s-1}J_0(P) \approx Z^{2s-1}U_0H^{\alpha/2}P^{\tau/2} = Z^{2s-1}PZS_{s-1}(P)H^{\alpha/2}P^{\tau/2}.$$

And

$$\lambda_s = \frac{\lambda_{s-1}}{(1+\theta)} + \frac{(2s-k\alpha/2)\theta}{(1+\theta)} + \frac{(1+\alpha/2+\tau/2)}{(1+\theta)}.$$

Or,

$$\Delta(s)=\frac{\Delta(s-1)}{(1+\theta)}+\frac{(1-\alpha/2)(k\theta-1)}{(1+\theta)}+\frac{\tau}{2(1+\theta)}. \tag{2.10}$$

By (2.10) and (2.9), there is

$$\Delta(s)=\Delta(s-1)\left(1-\frac{3\omega+(2-\beta)\delta}{2(k(\beta+1)+\delta+\omega)}\right)-\frac{k(\rho-\delta)}{2(k(\beta+1)+\delta+\omega)}$$

where $1-\alpha=\dfrac{\beta\cdot\Delta(s-1)}{k}, 1-\alpha-\tau=\dfrac{\omega\cdot\Delta(s-1)}{k}, \delta=\dfrac{\omega}{\beta}\left(1+\overline{\mu}\zeta\dfrac{\beta+1}{2}\right)\left(\dfrac{\beta-1}{2\beta}\right)^{k-1}$,

$\rho=1+\beta-\omega$.

And then,

$$\Delta(d+i)=\left(\Delta(d)+\frac{k(\rho-\delta)}{3\omega+(2-\beta)\delta}\right)\left(1-\frac{3\omega+(2-\beta)\delta}{2(k(\beta+1)+\delta+\omega)}\right)^{i}-\frac{k(\rho-\delta)}{3\omega+(2-\beta)\delta}. \tag{2.11}$$

Let

$$\Im(\beta,\omega)=\log\left(\frac{k(\rho-\delta)}{\Delta(d)(3\omega+(2-\beta)\delta)+k(\rho-\delta)}\right)\Big/\log\left(1-\frac{3\omega+(2-\beta)\delta}{2(k(\beta+1)+\delta+\omega)}\right). \tag{2.12}$$

We take $\beta\doteq 1/3$, and $a\doteq(-1)$. With (2.6) and (2.8), it follows

$$\theta_k>0\Rightarrow 1-\alpha>\tau_k\Rightarrow 1-\alpha-\tau>-\mu\sigma\frac{1-(2\sigma)^{k-1}}{1-2\sigma}\Rightarrow 1>-\overline{\mu}\sigma\frac{1-(2\sigma)^{k-1}}{1-2\sigma}.$$

And it is clear that, for larger $k$, and $a\doteq(-1)$, then

$\dfrac{1-(\sigma/a)^{k-1}}{1-(\sigma/a)}\Big/\dfrac{1-(2\sigma)^{k-1}}{1-2\sigma}=o(1)$, and $\zeta\doteq 2\sigma\dfrac{1-(2\sigma)^{k-1}}{(1-2\sigma)}$.

Hence, if $k$ is odd.

$$\delta=\frac{\omega}{\beta}\left(1+\overline{\mu}\zeta\frac{\beta+1}{2}\right)a^{k-1}\Rightarrow\delta\sim\frac{\omega}{\beta}\left(1-2\frac{\beta+1}{2}\right)a^{k-1}\sim-\omega a^{k-1}\sim-\omega,$$

If $k$ is even, we choose $\sigma$ and $\mu$ such that $\overline{\mu}\zeta\doteq-1$, equivalently, $\theta_k\doteq\dfrac{1-\alpha-\tau}{2k(1-\alpha)}$. Then

$$\delta=\frac{\omega}{\beta}\left(1+\overline{\mu}\zeta\frac{\beta+1}{2}\right)a^{k-1}\sim\frac{\omega}{\beta}\left(1-\frac{\beta+1}{2}\right)a^{k-1}\sim\omega a^{k-1}\sim-\omega. \qquad (\because\ \beta=1/3)$$

For the observation above, we may take $\sigma\doteq 1/2,\ \mu\doteq((-1)^k-3)(1-\alpha-\tau)/(2k)$.

From (2.12) we can know that, for larger $k$, there is

$$\Im(\beta,\omega) \sim \frac{2k(\beta+1)}{3\omega+(2-\beta)\delta}\log\left(1+\frac{3\omega+(2-\beta)\delta}{(\rho-\delta)}\frac{\Delta(d)}{k}\right) \sim \frac{2k}{\omega}\log\left(1+\omega\frac{\Delta(d)}{k}\right).$$

And then take $\Delta(d)$ such that $\Delta(d) \le 3k/7$, as in [1], it is easy to know $d \le \frac{4k}{3}\log\left(1+\frac{3}{4}\right)$,

and in total,

$$d+\Im(\beta,\omega) \le \frac{4k}{3}\log\left(\frac{7}{4}\right)+6k\log\left(\frac{8}{7}\right)=1.547343\cdot k. \tag{2.13}$$

Take $u=\Im(\beta,\omega)+d$, by Lemma 2

$$G(k) \le 2u+1. \tag{2.14}$$

Theorem 1 is followed.

The Proof of Theorem 2:

The proof is a routine as in [1], with computer, $\Delta(s)$ may be followed by the recursions (2.9) and (2.10), choosing optimally the parameters $\alpha$, $\tau$, $\mu$ and $\sigma$ in turn. As the restriction of the computation ability of PC, we have to take use of only four digits; the final results of $\Delta(s)$ are included in List 2.1, and the intermediate results of recursions are posted in the end of the paper as an appendix.

| $k$ | $s$ | $\Delta(s)$ | $k$ | $s$ | $\Delta(s)$ |
|---|---|---|---|---|---|
| 5 | 8 | 0.000000 | 13 | 20 | 0.000000 |
| 6 | 9 | 0.000000 | 14 | 22 | 0.000000 |
| 7 | 11 | 0.000000 | 15 | 23 | 0.000000 |
| 8 | 12 | 0.006974 | 16 | 25 | 0.000000 |
| 9 | 14 | 0.000000 | 17 | 26 | 0.000000 |
| 10 | 16 | 0.000000 | 18 | 28 | 0.000000 |
| 11 | 17 | 0.000000 | 19 | 30 | 0.000001 |
| 12 | 19 | 0.000000 | 20 | 31 | 0.000000 |

List 2.1

By Lemma.2, we take $v(k)$ as in the following list 2.2

| $k$ | 5 | 6 | 7 | 8 | 9 | 10 | 11 | 12 | 13 | 14 | 15 | 16 | 17 | 18 | 19 | 20 |
|---|---|---|---|---|---|---|---|---|---|---|---|---|---|---|---|---|
| $v(k)$ | 1 | 1 | 1 | 8 | 1 | 1 | 1 | 1 | 1 | 1 | 1 | 14 | 1 | 1 | 1 | 1 |

List 2.2

and the proof of Theorem 2 is finished.

As we seen, in this paper our main result is that:

If $k$ is larger, then there is an integer $u$ , $u \le 1.547343 \cdot k + o(k)$ , such that $\Delta(u) = 0$ .

Actually for the result there is more explicit and effective implication yet.

Let $P = N^{1/k}$ , $\mathscr{C} \subseteq [0, P]$ is an integer subset constructed as in [1], $|\mathscr{C}| = \rho \ge P^{1-\varepsilon}$, $\varepsilon$ is an arbitrary small positive real number.

For a natural number $n$ , denote by $r_n$ as the number of solutions of the following equation

$$x_1^k + \cdots + x_u^k = n, \quad x_i \in \mathscr{C},\ 1 \le i \le u. \tag{2.15}$$

Then clearly

$$\sum_n r_n = \rho^u . \tag{2.16}$$

And by Cauchy's inequality,

$$\sum_{r_d > 0} 1 \sum_n r_n^2 \ge \rho^{2u}. \tag{2.17}$$

Let $R = \sum_{r_d > 0} 1$, $\rho^{\lambda_u} = \sum_n r_n^2$ , then (2.17) may be rewritten as

$$R \cdot \rho^{\lambda_u} \ge \rho^{2u} ,$$

or

$$\log_\rho R + \lambda_u - 2u \ge 0 . \tag{2.18}$$

From (2.18), we can know that $\Delta(u)$ is actually a measure of the deviation of the inequality (2.18) underling the assuming of $R = \rho^k$. However, clearly the upper bound of $R$ is about $uN$, rather than $\rho^k (\le N)$ . Hence, we modify the definition of $\Delta(u)$ as

$$\overline{\Delta}(u) = \lambda_u - 2u + \log_\rho (uN) = \Delta(u) + o(k) + \log u / \log \rho. \tag{2.19}$$

For $P$ is great and term $\log u / \log \rho$ is minor, so our conclusion is also valid for $\overline{\Delta}(u)$ . If we define $\overline{G}(k)$ to be the least integer $u$ with $\overline{\Delta}(u) = 0$ , then it has

$$\overline{G}(k) \le 1.547343 \cdot k + o(k)\,. \qquad (2.20)$$

By the definition of $\overline{\Delta}(u)$, we can know that $\overline{\Delta}(u) = 0$ implies that $r_n = O(P^{u-k})$ for at least $O(N)$ ones of $n\ (\le N)$, we guess, for most ones of $n$ likely, except $k$ is a power of $2$.

So, it is natural to think that $\overline{G}(k)$ is closed to, even equal to $G(k)$, and we guess that if $k$ is not a power of $2$,

$$G(k) \le \overline{G}(k)\,. \qquad (2.21)$$

We have taken some concrete computations for $G(k)$ by PC.

For two positive integers $a, b,\ a \le b.$ write $X = [a, b]$, define $G(k; X)$ the least integer $s$ such that each integer in $X$ can be represented a sum of at most $s$ $k$-th powers of positive integers.

With the restriction of computing ability of PC, we just take $X = [1.6d \times 10^9, 3.2d \times 10^9]$, $1 \le d \le 113$, respectively. Denote $X_d = [1.6d \times 10^9, 3.2d \times 10^9]$, $1 \le d \le 113$, the results are as following

| $d$ | $G(5; X_d)$ | $G(6; X_d)$ | $G(7; X_d)$ |
|---|---|---|---|
| 1 | 9 | 15 | 20 |
| 2 | . | . | 19 |
| 3 | . | . | 19 |
| 4 | . | . | 18 |
| 5 | . | . | . |
| 6 | . | . | . |
| 7 | 9 | 14 | 18 |
| 8 | 9 | 14 | 18 |
| 9 | 9 | 14 | 17 |
| . | . | . | . |
| . | . | . | . |
| 28 | 9 | 14 | 16 |
| . | . | . | . |
| . | . | . | . |
| 43 | 8 | 14 | 16 |
| . | . | . | . |
| . | . | . | . |
| 113 | 8 | 13 | 15 |

List 2.3

By the data above, there is a preliminary estimate that $G(5) \le 8$, $G(6) \le 13$, $G(7) \le 15$, respectively. Especially, the first one is already very close and indicates $G(5) = 8, 7$, or $6$. For the rest two ones, the ranges will be further lessened yet as the size of $X$ increase, and we guess that $G(k; X)$ will be near to $G(k)$ when $X = [10^{2k}, 2\times 10^{2k}]$, or greater.

In fact, the calculations of $G(k; X)$ will be helpful and meaningful not only to estimate the upper bound of $G(k)$, but also to determine the exact value of $G(k)$.

In the process of computing of $G(k; X)$, it will be encountered the case

$$G(k; X_d) = G(k; X_{d+1}) = \cdots = G(k, X_{d+r}) = \cdots (= v, \text{ say}). \tag{2.22}$$

e.g. $G(5; X_{43}) = \cdots = G(5; X_{113}) = 8$ in List 2.3.

At this time, it is natural to ask whether the value $v$ is already equal to $G(k)$, or will be decreased further.

Let $k$, $s$ and $X$ as before, and denote $X'(s)$ as the subset of set $X$ in which the elements can be represented as a sum of at most $s$ $k$-th powers of the natural numbers, and $X''(s) = X \setminus X'(s)$. In recursion calculation of $G(k; X)$, we are able to calculate the size $|X''(s)|$ for each $s$. Obviously, $|X''(s)|$ will be decreased or increased as the size of $X$ increase.

In respect to (2.22), clearly, there are

$$X''_d(v) = \cdots = X''_{d+r}(v) = \varnothing. \tag{2.23}$$

So, the trend of the sizes of the sets $X''_d(v-1), \cdots, X''_{d+r}(v-1)$ will play an important role for the next value of $G(k; X)$. Say in detail, that decreasing will imply $G(k) \le v-1$, and otherwise, increasing (or undulating ) will imply $G(k) = v$. In the following is the data of $G(5; X_d)$

| $d$ | 43 | 44 | 45 | ······ | 113 |
|---|---|---|---|---|---|
| $\lvert X''_d(7) \rvert$ | 6470475034 | 6524570166 | 6576538986 | ······ | 8482079778 |

List 2.4

Recently, we have also taken the test for $G(3)$. Here $X_d$ has been replaced with $Y_d$, $Y_d =$

$[3.2d \times 10^7, 6.4d \times 10^7]$, the test results are as following

| $d$ | $G(3;Y_d)$ | $\lvert Y_d''(4) \rvert$ |
|---|---|---|
| 1 | 5 | 2540807 |
| 2 | . | 3870200 |
| . | . | . |
| . | . | . |
| . | . | . |
| 109 | 5 | 14769036 |
| 110 | 5 | 14769899 |
| 111 | 5 | 14769703 |
| . | . | . |
| . | . | . |
| . | . | . |
| 1000 | 5 | 7397806 |
| . | . | . |
| . | . | . |
| . | . | . |
| 2000 | 5 | 4062610 |

List 2.5

From List 2.5, we have seen that the sizes of $\lvert Y_d''(4) \rvert$ are increasing from $d = 1$ to $d = 110$, and then decreasing to $d = 2000$, that is, which are undulated , rather than monotone that I thought before. (the case of $G(5; X_d)$ maybe is also similar, that is, undulant, only we are not able to do more).

Although we don't know what will happen next, we prefer that

$$G(5) = 8, \quad G(3) = 5. \tag{2.24}$$

Besides, we noticed that the values of $G(k; X)$ presented in List 2.3 and estimation (2.20) for $\overline{G}(k)$, and (2.24) as well all seem to approach $\varsigma(k)$ defined in Lemma 2, that make us guess that

$$G(k), \overline{G}(k) \le \varsigma(k) + 1. \tag{2.25}$$

It will be obliged for the universities, institutes and companies with the modern computers such Google, IBM, etc. to present further results for $G(k; X)$. A reference code is posted in the

Appendix II.

## Appendix The Intermediate Results of Recursions for Theorem 2

| $k$ | $s$ | $\alpha$ | $\tau$ | $\mu$ | $\sigma$ | $\theta$ | $\Delta(s)$ |
|---|---|---|---|---|---|---|---|
| 5 | 3 | 0.000000 | 0.000000 | 0.000000 | 0.000000 | 0.125120 | 2.333618 |
| | 4 | 0.000000 | 0.000000 | 0.000000 | 0.000000 | 0.136680 | 1.774482 |
| | 5 | 0.880000 | 0.000000 | -0.056800 | 0.500300 | 0.000002 | 1.214485 |
| | 6 | 0.918000 | 0.000000 | -0.038300 | 0.502600 | 0.000005 | 0.673494 |
| | 7 | 0.955000 | 0.000000 | -0.020800 | 0.501400 | 0.000010 | 0.151018 |
| | 8 | -0.263000 | 0.341000 | -0.341000 | 0.500000 | 0.143170 | 0.000000 |
| 6 | 3 | 0.000000 | 0.000000 | 0.000000 | 0.000000 | 0.100009 | 3.272749 |
| | 4 | 0.000000 | 0.000000 | 0.000000 | 0.000000 | 0.107879 | 2.635686 |
| | 5 | 0.000000 | 0.000000 | 0.000000 | 0.000000 | 0.115887 | 2.088927 |
| | 6 | 0.883000 | 0.000000 | -0.021500 | 0.497600 | 0.000554 | 1.531434 |
| | 7 | 0.914000 | 0.000000 | -0.015233 | 0.502500 | 0.000889 | 0.990449 |
| | 8 | 0.944000 | 0.000000 | -0.010133 | 0.495600 | 0.002130 | 0.468199 |
| | 9 | 0.438000 | 0.518000 | -0.518000 | 0.000000 | 0.011481 | 0.000000 |
| 7 | 3 | 0.000000 | 0.000000 | 0.000000 | 0.000000 | 0.083334 | 4.230771 |
| | 4 | 0.000000 | 0.000000 | 0.000000 | 0.000000 | 0.089044 | 3.538957 |
| | 5 | 0.000000 | 0.000000 | 0.000000 | 0.000000 | 0.094897 | 2.925605 |
| | 6 | 0.860000 | 0.000000 | -0.043300 | 0.504700 | 0.000172 | 2.355885 |
| | 7 | 0.887000 | 0.000000 | -0.035286 | 0.503800 | 0.000423 | 1.800272 |
| | 8 | 0.914000 | 0.000000 | -0.026471 | 0.505000 | 0.000009 | 1.257295 |
| | 9 | 0.940000 | 0.000000 | -0.019143 | 0.499300 | 0.000004 | 0.727308 |
| | 10 | 0.965000 | 0.000000 | -0.011000 | 0.503600 | 0.000774 | 0.212447 |
| | 11 | -0.462000 | 0.578000 | -4.046000 | 0.100000 | 0.084664 | 0.000000 |
| 8 | 3 | 0.000000 | 0.000000 | 0.000000 | 0.000000 | 0.071429 | 5.200000 |
| | 4 | 0.000000 | 0.000000 | 0.000000 | 0.000000 | 0.075758 | 4.467607 |
| | 5 | 0.000000 | 0.000000 | 0.000000 | 0.000000 | 0.080209 | 3.804151 |
| | 6 | 0.000000 | 0.000000 | 0.000000 | 0.000000 | 0.084719 | 3.209961 |
| | 7 | 0.866000 | 0.000000 | -0.017650 | 0.503200 | 0.000114 | 2.643177 |
| | 8 | 0.889000 | 0.000000 | -0.015175 | 0.495700 | 0.001141 | 2.090363 |
| | 9 | 0.912000 | 0.000000 | -0.011300 | 0.502800 | 0.001596 | 1.550834 |
| | 10 | 0.935000 | 0.000000 | -0.008625 | 0.500400 | 0.000778 | 1.020854 |
| | 11 | 0.957000 | 0.000000 | -0.005775 | 0.496800 | 0.001729 | 0.505694 |
| | 12 | 0.978000 | 0.021000 | -0.000125 | 0.498700 | 0.000436 | 0.006974 |
| 9 | 3 | 0.000000 | 0.000000 | 0.000000 | 0.000000 | 0.062500 | 6.176471 |
| | 4 | 0.000000 | 0.000000 | 0.000000 | 0.000000 | 0.065891 | 5.412834 |
| | 5 | 0.000000 | 0.000000 | 0.000000 | 0.000000 | 0.069383 | 4.710455 |
| | 6 | 0.000000 | 0.000000 | 0.000000 | 0.000000 | 0.072937 | 4.070034 |
| | 7 | 0.000000 | 0.000000 | 0.000000 | 0.000000 | 0.076512 | 3.491500 |
| | 8 | 0.870000 | 0.000000 | -0.030689 | 0.503100 | 0.000616 | 2.927829 |
| | 9 | 0.891000 | 0.000000 | -0.026822 | 0.498300 | 0.000601 | 2.374900 |
| | 10 | 0.912000 | 0.000000 | -0.021956 | 0.495700 | 0.000007 | 1.830920 |
| | 11 | 0.932000 | 0.000000 | -0.015811 | 0.504400 | 0.000204 | 1.297635 |
| | 12 | 0.951000 | 0.000000 | -0.012389 | 0.496900 | 0.003570 | 0.787178 |
| | 13 | 0.970000 | 0.027000 | -0.000767 | 0.496500 | 0.000562 | 0.288121 |
| | 14 | -0.979000 | 1.099000 | -2.198000 | 0.200000 | 0.048628 | 0.000000 |
| 10 | 3 | 0.000000 | 0.000000 | 0.000000 | 0.000000 | 0.055556 | 7.157895 |
| | 4 | 0.000000 | 0.000000 | 0.000000 | 0.000000 | 0.058282 | 6.369489 |
| | 5 | 0.000000 | 0.000000 | 0.000000 | 0.000000 | 0.061089 | 5.636078 |
| | 6 | 0.000000 | 0.000000 | 0.000000 | 0.000000 | 0.063955 | 4.958505 |
| | 7 | 0.000000 | 0.000000 | 0.000000 | 0.000000 | 0.066852 | 4.337081 |
| | 8 | 0.000000 | 0.000000 | 0.000000 | 0.000000 | 0.069750 | 3.771515 |
| | 9 | 0.874000 | 0.000000 | -0.013900 | 0.495800 | 0.000287 | 3.209210 |
| | 10 | 0.893000 | 0.000000 | -0.011400 | 0.500400 | 0.000009 | 2.655736 |
| | 11 | 0.911000 | 0.000000 | -0.009500 | 0.498200 | 0.000978 | 2.114493 |
| | 12 | 0.929000 | 0.000000 | -0.007700 | 0.495900 | 0.001400 | 1.584271 |
| | 13 | 0.947000 | 0.000000 | -0.005800 | 0.496200 | 0.000595 | 1.060272 |
| | 14 | 0.964000 | 0.035000 | -0.000100 | 0.502300 | 0.000122 | 0.560337 |
| | 15 | 0.981000 | 0.018000 | -0.000100 | 0.502500 | 0.000213 | 0.060908 |
| | 16 | -0.016000 | 0.113000 | -0.452000 | 0.400000 | 0.088352 | 0.000000 |
| 11 | 3 | 0.000000 | 0.000000 | 0.000000 | 0.000000 | 0.050000 | 8.142857 |
| | 4 | 0.000000 | 0.000000 | 0.000000 | 0.000000 | 0.052239 | 7.334347 |
| | 5 | 0.000000 | 0.000000 | 0.000000 | 0.000000 | 0.054542 | 6.575661 |
| | 6 | 0.000000 | 0.000000 | 0.000000 | 0.000000 | 0.056897 | 5.867674 |
| | 7 | 0.000000 | 0.000000 | 0.000000 | 0.000000 | 0.059285 | 5.210882 |
| | 8 | 0.000000 | 0.000000 | 0.000000 | 0.000000 | 0.061687 | 4.605349 |
| | 9 | 0.860000 | 0.000000 | -0.027455 | 0.499400 | 0.000471 | 4.036402 |
| | 10 | 0.877000 | 0.000000 | -0.025264 | 0.495400 | 0.000950 | 3.477465 |
| | 11 | 0.894000 | 0.000000 | -0.020873 | 0.499400 | 0.001014 | 2.927664 |
| | 12 | 0.911000 | 0.000000 | -0.017982 | 0.496600 | 0.000475 | 2.384876 |
| | 13 | 0.927000 | 0.000000 | -0.014773 | 0.497300 | 0.001953 | 1.856275 |
| | 14 | 0.943000 | 0.000000 | -0.011764 | 0.495700 | 0.002733 | 1.340002 |
| | 15 | 0.959000 | 0.000000 | -0.008155 | 0.498900 | 0.001844 | 0.828531 |
| | 16 | 0.974000 | 0.025000 | -0.000182 | 0.505000 | 0.000385 | 0.330075 |
| | 17 | -0.754000 | 1.129000 | -1.129000 | 0.200000 | 0.031850 | 0.000000 |
| 12 | 3 | 0.000000 | 0.000000 | 0.000000 | 0.000000 | 0.045455 | 9.130435 |
| | 4 | 0.000000 | 0.000000 | 0.000000 | 0.000000 | 0.047325 | 8.305287 |
| | 5 | 0.000000 | 0.000000 | 0.000000 | 0.000000 | 0.049248 | 7.525642 |
| | 6 | 0.000000 | 0.000000 | 0.000000 | 0.000000 | 0.051215 | 6.792350 |
| | 7 | 0.000000 | 0.000000 | 0.000000 | 0.000000 | 0.053213 | 6.105989 |
| | 8 | 0.000000 | 0.000000 | 0.000000 | 0.000000 | 0.055230 | 5.466819 |
| | 9 | 0.000000 | 0.000000 | 0.000000 | 0.000000 | 0.057251 | 4.874750 |
| | 10 | 0.864000 | 0.000000 | -0.012033 | 0.498300 | 0.000925 | 4.309069 |
| | 11 | 0.880000 | 0.000000 | -0.010300 | 0.501400 | 0.000489 | 3.750520 |
| | 12 | 0.895000 | 0.000000 | -0.009250 | 0.497800 | 0.001785 | 3.204134 |
| | 13 | 0.910000 | 0.000000 | -0.008100 | 0.495500 | 0.002725 | 2.669679 |
| | 14 | 0.925000 | 0.000000 | -0.006450 | 0.499400 | 0.002776 | 2.144131 |
| | 15 | 0.940000 | 0.000000 | -0.005400 | 0.496000 | 0.001642 | 1.621911 |
| | 16 | 0.954000 | 0.045000 | -0.000083 | 0.500800 | 0.000115 | 1.122003 |
| | 17 | 0.968000 | 0.031000 | -0.000083 | 0.500000 | 0.000206 | 0.622652 |
| | 18 | 0.982000 | 0.017000 | -0.000083 | 0.497900 | 0.000523 | 0.125282 |
| | 19 | -0.240000 | 0.293000 | -2.344000 | 0.100000 | 0.063111 | 0.000000 |
| 13 | 3 | 0.000000 | 0.000000 | 0.000000 | 0.000000 | 0.041667 | 10.120000 |
| | 4 | 0.000000 | 0.000000 | 0.000000 | 0.000000 | 0.043253 | 9.280862 |
| | 5 | 0.000000 | 0.000000 | 0.000000 | 0.000000 | 0.044882 | 8.483567 |
| | 6 | 0.000000 | 0.000000 | 0.000000 | 0.000000 | 0.046547 | 7.728921 |
| | 7 | 0.000000 | 0.000000 | 0.000000 | 0.000000 | 0.048242 | 7.017526 |
| | 8 | 0.000000 | 0.000000 | 0.000000 | 0.000000 | 0.049956 | 6.349748 |
| | 9 | 0.000000 | 0.000000 | 0.000000 | 0.000000 | 0.051680 | 5.725686 |
| | 10 | 0.000000 | 0.000000 | 0.000000 | 0.000000 | 0.053403 | 5.145155 |
| | 11 | 0.868000 | 0.000000 | -0.022808 | 0.495100 | 0.000068 | 4.579344 |
| | 12 | 0.882000 | 0.000000 | -0.018754 | 0.502300 | 0.000990 | 4.023554 |
| | 13 | 0.896000 | 0.000000 | -0.018200 | 0.495000 | 0.001718 | 3.477907 |
| | 14 | 0.910000 | 0.000000 | -0.015646 | 0.495600 | 0.002021 | 2.941283 |
| | 15 | 0.924000 | 0.000000 | -0.012592 | 0.499300 | 0.001641 | 2.410802 |
| | 16 | 0.938000 | 0.000000 | -0.010038 | 0.500600 | 0.000559 | 1.882609 |
| | 17 | 0.951000 | 0.000000 | -0.008438 | 0.496700 | 0.003493 | 1.377114 |
| | 18 | 0.964000 | 0.035000 | -0.000154 | 0.505000 | 0.000142 | 0.877448 |
| | 19 | 0.977000 | 0.022000 | -0.000154 | 0.505000 | 0.000253 | 0.378536 |
| | 20 | -0.030000 | 0.765000 | -1.530000 | 0.100000 | 0.019247 | 0.000000 |
| 14 | 3 | 0.000000 | 0.000000 | 0.000000 | 0.000000 | 0.038462 | 11.111111 |
| | 4 | 0.000000 | 0.000000 | 0.000000 | 0.000000 | 0.039823 | 10.260047 |
| | 5 | 0.000000 | 0.000000 | 0.000000 | 0.000000 | 0.041220 | 9.447694 |
| | 6 | 0.000000 | 0.000000 | 0.000000 | 0.000000 | 0.042648 | 8.674803 |
| | 7 | 0.000000 | 0.000000 | 0.000000 | 0.000000 | 0.044102 | 7.941973 |
| | 8 | 0.000000 | 0.000000 | 0.000000 | 0.000000 | 0.045575 | 7.249620 |
| | 9 | 0.000000 | 0.000000 | 0.000000 | 0.000000 | 0.047060 | 6.597958 |
| | 10 | 0.000000 | 0.000000 | 0.000000 | 0.000000 | 0.048549 | 5.986978 |
| | 11 | 0.857000 | 0.000000 | -0.011014 | 0.497100 | 0.000726 | 5.417353 |
| | 12 | 0.871000 | 0.000000 | -0.009614 | 0.500400 | 0.000016 | 4.852901 |
| | 13 | 0.884000 | 0.000000 | -0.009086 | 0.495700 | 0.000925 | 4.298151 |
| | 14 | 0.897000 | 0.000000 | -0.008057 | 0.495200 | 0.001562 | 3.752850 |
| | 15 | 0.910000 | 0.000000 | -0.006829 | 0.497400 | 0.001800 | 3.215793 |
| | 16 | 0.923000 | 0.000000 | -0.006000 | 0.495700 | 0.001369 | 2.683939 |
| | 17 | 0.936000 | 0.000000 | -0.004571 | 0.503000 | 0.000298 | 2.153515 |
| | 18 | 0.948000 | 0.000000 | -0.004014 | 0.495300 | 0.003702 | 1.648671 |
| | 19 | 0.960000 | 0.039000 | -0.000071 | 0.500500 | 0.000123 | 1.148927 |
| | 20 | 0.972000 | 0.027000 | -0.000071 | 0.499900 | 0.000213 | 0.649823 |
| | 21 | 0.984000 | 0.015000 | -0.000071 | 0.498600 | 0.000507 | 0.152854 |
| | 22 | 0.112000 | 0.257000 | -0.771000 | 0.100000 | 0.050140 | 0.000000 |
| 15 | 3 | 0.000000 | 0.000000 | 0.000000 | 0.000000 | 0.035714 | 12.103448 |
| | 4 | 0.000000 | 0.000000 | 0.000000 | 0.000000 | 0.036896 | 11.242099 |
| | 5 | 0.000000 | 0.000000 | 0.000000 | 0.000000 | 0.038107 | 10.416751 |
| | 6 | 0.000000 | 0.000000 | 0.000000 | 0.000000 | 0.039344 | 9.628104 |
| | 7 | 0.000000 | 0.000000 | 0.000000 | 0.000000 | 0.040604 | 8.876733 |
| | 8 | 0.000000 | 0.000000 | 0.000000 | 0.000000 | 0.041882 | 8.163076 |
| | 9 | 0.000000 | 0.000000 | 0.000000 | 0.000000 | 0.043172 | 7.487410 |
| | 10 | 0.000000 | 0.000000 | 0.000000 | 0.000000 | 0.044469 | 6.849842 |
| | 11 | 0.000000 | 0.000000 | 0.000000 | 0.000000 | 0.045767 | 6.250290 |
| | 12 | 0.861000 | 0.000000 | -0.020633 | 0.495900 | 0.000129 | 5.681159 |
| | 13 | 0.873000 | 0.000000 | -0.019233 | 0.495100 | 0.001341 | 5.122124 |
| | 14 | 0.886000 | 0.000000 | -0.017100 | 0.495300 | 0.000286 | 4.566206 |
| | 15 | 0.898000 | 0.000000 | -0.014900 | 0.497500 | 0.001170 | 4.020171 |
| | 16 | 0.910000 | 0.000000 | -0.013100 | 0.497900 | 0.001744 | 3.483353 |
| | 17 | 0.922000 | 0.000000 | -0.010800 | 0.501300 | 0.001822 | 2.953704 |
| | 18 | 0.934000 | 0.000000 | -0.009800 | 0.496400 | 0.001250 | 2.427663 |
| | 19 | 0.946000 | 0.000000 | -0.007800 | 0.497800 | 0.000177 | 1.901726 |
| | 20 | 0.957000 | 0.042000 | -0.000133 | 0.504300 | 0.000100 | 1.401870 |
| | 21 | 0.968000 | 0.031000 | -0.000133 | 0.504500 | 0.000204 | 0.902767 |
| | 22 | 0.979000 | 0.020000 | -0.000133 | 0.504500 | 0.000484 | 0.405774 |
| | 23 | -0.940000 | 1.559000 | -1.559000 | 0.100000 | 0.012913 | 0.000000 |
| 16 | 3 | 0.000000 | 0.000000 | 0.000000 | 0.000000 | 0.033333 | 13.096774 |
| | 4 | 0.000000 | 0.000000 | 0.000000 | 0.000000 | 0.034368 | 12.226463 |
| | 5 | 0.000000 | 0.000000 | 0.000000 | 0.000000 | 0.035428 | 11.389793 |
| | 6 | 0.000000 | 0.000000 | 0.000000 | 0.000000 | 0.036510 | 10.587406 |
| | 7 | 0.000000 | 0.000000 | 0.000000 | 0.000000 | 0.037612 | 9.819853 |
| | 8 | 0.000000 | 0.000000 | 0.000000 | 0.000000 | 0.038730 | 9.087570 |
| | 9 | 0.000000 | 0.000000 | 0.000000 | 0.000000 | 0.039860 | 8.390873 |
| | 10 | 0.000000 | 0.000000 | 0.000000 | 0.000000 | 0.040999 | 7.729937 |
| | 11 | 0.000000 | 0.000000 | 0.000000 | 0.000000 | 0.042141 | 7.104789 |
| | 12 | 0.000000 | 0.000000 | 0.000000 | 0.000000 | 0.043281 | 6.515297 |
| | 13 | 0.864000 | 0.000000 | -0.009200 | 0.497100 | 0.000457 | 5.948732 |
| | 14 | 0.876000 | 0.000000 | -0.007850 | 0.501700 | 0.000112 | 5.387136 |
| | 15 | 0.887000 | 0.000000 | -0.007662 | 0.496100 | 0.001801 | 4.837958 |
| | 16 | 0.899000 | 0.000000 | -0.007012 | 0.495400 | 0.000493 | 4.289686 |
| | 17 | 0.910000 | 0.000000 | -0.006025 | 0.496900 | 0.001864 | 3.753943 |
| | 18 | 0.921000 | 0.000000 | -0.005337 | 0.496000 | 0.002735 | 3.229219 |
| | 19 | 0.932000 | 0.000000 | -0.004450 | 0.498000 | 0.002911 | 2.712198 |
| | 20 | 0.943000 | 0.000000 | -0.003663 | 0.499400 | 0.002371 | 2.198531 |
| | 21 | 0.954000 | 0.000000 | -0.003075 | 0.497400 | 0.001092 | 1.682835 |
| | 22 | 0.964000 | 0.035000 | -0.000063 | 0.499300 | 0.000185 | 1.183647 |
| | 23 | 0.975000 | 0.024000 | -0.000063 | 0.500700 | 0.000149 | 0.684270 |
| | 24 | 0.985000 | 0.014000 | -0.000063 | 0.496500 | 0.000708 | 0.189388 |
| | 25 | -0.380000 | 0.505000 | -0.505000 | 0.000000 | 0.039292 | 0.000000 |
| 17 | 3 | 0.000000 | 0.000000 | 0.000000 | 0.000000 | 0.031250 | 14.090909 |
| | 4 | 0.000000 | 0.000000 | 0.000000 | 0.000000 | 0.032164 | 13.212722 |
| | 5 | 0.000000 | 0.000000 | 0.000000 | 0.000000 | 0.033099 | 12.366098 |
| | 6 | 0.000000 | 0.000000 | 0.000000 | 0.000000 | 0.034053 | 11.551631 |
| | 7 | 0.000000 | 0.000000 | 0.000000 | 0.000000 | 0.035024 | 10.769837 |
| | 8 | 0.000000 | 0.000000 | 0.000000 | 0.000000 | 0.036010 | 10.021148 |
| | 9 | 0.000000 | 0.000000 | 0.000000 | 0.000000 | 0.037008 | 9.305892 |
| | 10 | 0.000000 | 0.000000 | 0.000000 | 0.000000 | 0.038014 | 8.624289 |
| | 11 | 0.000000 | 0.000000 | 0.000000 | 0.000000 | 0.039025 | 7.976438 |
| | 12 | 0.000000 | 0.000000 | 0.000000 | 0.000000 | 0.040038 | 7.362309 |
| | 13 | 0.000000 | 0.000000 | 0.000000 | 0.000000 | 0.041047 | 6.781738 |
| | 14 | 0.867000 | 0.000000 | -0.016347 | 0.499800 | 0.000040 | 6.215376 |
| | 15 | 0.878000 | 0.000000 | -0.015953 | 0.496200 | 0.000215 | 5.655212 |
| | 16 | 0.889000 | 0.000000 | -0.013559 | 0.500300 | 0.000208 | 5.100615 |
| | 17 | 0.899000 | 0.000000 | -0.013582 | 0.495100 | 0.002519 | 4.562197 |
| | 18 | 0.910000 | 0.000000 | -0.012088 | 0.495000 | 0.001501 | 4.025063 |
| | 19 | 0.921000 | 0.000000 | -0.010294 | 0.496400 | 0.000197 | 3.486681 |
| | 20 | 0.931000 | 0.000000 | -0.009018 | 0.496800 | 0.002367 | 2.966667 |
| | 21 | 0.941000 | 0.000000 | -0.007641 | 0.497600 | 0.003601 | 2.460723 |
| | 22 | 0.951000 | 0.048000 | -0.000118 | 0.503300 | 0.000081 | 1.960788 |
| | 23 | 0.961000 | 0.038000 | -0.000118 | 0.503300 | 0.000094 | 1.460982 |
| | 24 | 0.971000 | 0.028000 | -0.000118 | 0.503200 | 0.000110 | 0.961336 |
| | 25 | 0.981000 | 0.018000 | -0.000118 | 0.503000 | 0.000109 | 0.461726 |
| | 26 | -0.231000 | 1.132000 | -1.132000 | 0.000000 | 0.004629 | 0.000000 |
| 18 | 3 | 0.000000 | 0.000000 | 0.000000 | 0.000000 | 0.029412 | 15.085714 |
| | 4 | 0.000000 | 0.000000 | 0.000000 | 0.000000 | 0.030225 | 14.200551 |
| | 5 | 0.000000 | 0.000000 | 0.000000 | 0.000000 | 0.031055 | 13.345110 |
| | 6 | 0.000000 | 0.000000 | 0.000000 | 0.000000 | 0.031903 | 12.519940 |
| | 7 | 0.000000 | 0.000000 | 0.000000 | 0.000000 | 0.032765 | 11.725526 |
| | 8 | 0.000000 | 0.000000 | 0.000000 | 0.000000 | 0.033641 | 10.962283 |
| | 9 | 0.000000 | 0.000000 | 0.000000 | 0.000000 | 0.034528 | 10.230544 |
| | 10 | 0.000000 | 0.000000 | 0.000000 | 0.000000 | 0.035423 | 9.530554 |
| | 11 | 0.000000 | 0.000000 | 0.000000 | 0.000000 | 0.036323 | 8.862459 |
| | 12 | 0.000000 | 0.000000 | 0.000000 | 0.000000 | 0.037227 | 8.226302 |
| | 13 | 0.000000 | 0.000000 | 0.000000 | 0.000000 | 0.038130 | 7.622011 |
| | 14 | 0.858000 | 0.000000 | -0.008589 | 0.496200 | 0.001639 | 7.056290 |
| | 15 | 0.869000 | 0.000000 | -0.007778 | 0.497800 | 0.000652 | 6.493195 |
| | 16 | 0.879000 | 0.000000 | -0.007422 | 0.495400 | 0.001725 | 5.939850 |
| | 17 | 0.890000 | 0.000000 | -0.006611 | 0.497700 | 0.000021 | 5.384946 |
| | 18 | 0.900000 | 0.000000 | -0.006156 | 0.495700 | 0.000740 | 4.838690 |
| | 19 | 0.910000 | 0.000000 | -0.005400 | 0.496900 | 0.001186 | 4.300227 |
| | 20 | 0.920000 | 0.000000 | -0.004844 | 0.496300 | 0.001223 | 3.767504 |
| | 21 | 0.930000 | 0.000000 | -0.003989 | 0.500000 | 0.000880 | 3.238129 |
| | 22 | 0.940000 | 0.000000 | -0.003433 | 0.500300 | 0.000137 | 2.709063 |
| | 23 | 0.949000 | 0.050000 | -0.000056 | 0.500200 | 0.000087 | 2.209196 |
| | 24 | 0.959000 | 0.000000 | -0.002378 | 0.499200 | 0.000561 | 1.693003 |
| | 25 | 0.968000 | 0.031000 | -0.000056 | 0.499800 | 0.000167 | 1.193857 |
| | 26 | 0.977000 | 0.022000 | -0.000056 | 0.497900 | 0.000379 | 0.696585 |
| | 27 | 0.987000 | 0.012000 | -0.000056 | 0.501000 | 0.000169 | 0.197589 |
| | 28 | 0.455000 | 0.205000 | -0.205000 | 0.000000 | 0.033974 | 0.000000 |
| 19 | 3 | 0.000000 | 0.000000 | 0.000000 | 0.000000 | 0.027778 | 16.081081 |
| | 4 | 0.000000 | 0.000000 | 0.000000 | 0.000000 | 0.028505 | 15.189695 |
| | 5 | 0.000000 | 0.000000 | 0.000000 | 0.000000 | 0.029249 | 14.326392 |
| | 6 | 0.000000 | 0.000000 | 0.000000 | 0.000000 | 0.030006 | 13.491676 |
| | 7 | 0.000000 | 0.000000 | 0.000000 | 0.000000 | 0.030777 | 12.686002 |
| | 8 | 0.000000 | 0.000000 | 0.000000 | 0.000000 | 0.031560 | 11.909768 |
| | 9 | 0.000000 | 0.000000 | 0.000000 | 0.000000 | 0.032352 | 11.163303 |
| | 10 | 0.000000 | 0.000000 | 0.000000 | 0.000000 | 0.033153 | 10.446864 |
| | 11 | 0.000000 | 0.000000 | 0.000000 | 0.000000 | 0.033959 | 9.760628 |
| | 12 | 0.000000 | 0.000000 | 0.000000 | 0.000000 | 0.034770 | 9.104685 |
| | 13 | 0.000000 | 0.000000 | 0.000000 | 0.000000 | 0.035581 | 8.479035 |
| | 14 | 0.000000 | 0.000000 | 0.000000 | 0.000000 | 0.036391 | 7.883577 |
| | 15 | 0.861000 | 0.000000 | -0.016532 | 0.495800 | 0.001265 | 7.318508 |
| | 16 | 0.871000 | 0.000000 | -0.015479 | 0.495300 | 0.001186 | 6.758714 |
| | 17 | 0.881000 | 0.000000 | -0.014226 | 0.495400 | 0.000908 | 6.203235 |
| | 18 | 0.891000 | 0.000000 | -0.012874 | 0.495900 | 0.000372 | 5.650551 |
| | 19 | 0.900000 | 0.000000 | -0.011926 | 0.495800 | 0.002300 | 5.112825 |
| | 20 | 0.910000 | 0.000000 | -0.010774 | 0.495300 | 0.000855 | 4.572768 |
| | 21 | 0.919000 | 0.000000 | -0.009526 | 0.496600 | 0.002538 | 4.048059 |
| | 22 | 0.928000 | 0.000000 | -0.008479 | 0.496700 | 0.003560 | 3.535730 |
| | 23 | 0.937000 | 0.000000 | -0.007132 | 0.498900 | 0.003985 | 3.032388 |
| | 24 | 0.946000 | 0.053000 | -0.000105 | 0.502700 | 0.000072 | 2.532426 |
| | 25 | 0.955000 | 0.044000 | -0.000105 | 0.502600 | 0.000075 | 2.032521 |
| | 26 | 0.964000 | 0.035000 | -0.000105 | 0.502500 | 0.000071 | 1.532609 |
| | 27 | 0.973000 | 0.000000 | -0.003242 | 0.495200 | 0.001062 | 1.028379 |
| | 28 | 0.981000 | 0.018000 | -0.000105 | 0.502800 | 0.000533 | 0.532755 |
| | 29 | 0.990000 | 0.009000 | -0.000105 | 0.502800 | 0.001153 | 0.043264 |
| | 30 | -0.055000 | 0.087000 | -0.087000 | 0.000000 | 0.048187 | 0.000001 |
| 20 | 3 | 0.000000 | 0.000000 | 0.000000 | 0.000000 | 0.026316 | 17.076923 |
| | 4 | 0.000000 | 0.000000 | 0.000000 | 0.000000 | 0.026971 | 16.179953 |
| | 5 | 0.000000 | 0.000000 | 0.000000 | 0.000000 | 0.027640 | 15.309594 |
| | 6 | 0.000000 | 0.000000 | 0.000000 | 0.000000 | 0.028321 | 14.466313 |
| | 7 | 0.000000 | 0.000000 | 0.000000 | 0.000000 | 0.029014 | 13.650536 |
| | 8 | 0.000000 | 0.000000 | 0.000000 | 0.000000 | 0.029717 | 12.862638 |
| | 9 | 0.000000 | 0.000000 | 0.000000 | 0.000000 | 0.030430 | 12.102943 |
| | 10 | 0.000000 | 0.000000 | 0.000000 | 0.000000 | 0.031150 | 11.371713 |
| | 11 | 0.000000 | 0.000000 | 0.000000 | 0.000000 | 0.031876 | 10.669142 |
| | 12 | 0.000000 | 0.000000 | 0.000000 | 0.000000 | 0.032606 | 9.995354 |
| | 13 | 0.000000 | 0.000000 | 0.000000 | 0.000000 | 0.033338 | 9.350396 |
| | 14 | 0.000000 | 0.000000 | 0.000000 | 0.000000 | 0.034071 | 8.734233 |
| | 15 | 0.000000 | 0.000000 | 0.000000 | 0.000000 | 0.034802 | 8.146746 |
| | 16 | 0.864000 | 0.000000 | -0.007700 | 0.495000 | 0.000432 | 7.580378 |
| | 17 | 0.873000 | 0.000000 | -0.007050 | 0.495500 | 0.001458 | 7.023069 |
| | 18 | 0.882000 | 0.000000 | -0.006300 | 0.497300 | 0.002298 | 6.474881 |
| | 19 | 0.892000 | 0.000000 | -0.005900 | 0.497000 | 0.000196 | 5.921891 |
| | 20 | 0.901000 | 0.000000 | -0.005550 | 0.495300 | 0.000847 | 5.377146 |
| | 21 | 0.910000 | 0.000000 | -0.004900 | 0.496600 | 0.001168 | 4.839226 |
| | 22 | 0.919000 | 0.000000 | -0.004350 | 0.497300 | 0.001180 | 4.306403 |
| | 23 | 0.928000 | 0.000000 | -0.003800 | 0.498400 | 0.000878 | 3.776499 |
| | 24 | 0.937000 | 0.000000 | -0.003250 | 0.499900 | 0.000231 | 3.246707 |
| | 25 | 0.945000 | 0.054000 | -0.000050 | 0.500000 | 0.000080 | 2.746829 |
| | 26 | 0.954000 | 0.000000 | -0.002500 | 0.496700 | 0.001356 | 2.234984 |
| | 27 | 0.962000 | 0.037000 | -0.000050 | 0.499700 | 0.000137 | 1.735664 |
| | 28 | 0.971000 | 0.000000 | -0.001550 | 0.497900 | 0.000711 | 1.227609 |
| | 29 | 0.979000 | 0.020000 | -0.000050 | 0.499200 | 0.000303 | 0.729986 |
| | 30 | 0.987000 | 0.012000 | -0.000050 | 0.495400 | 0.000844 | 0.237836 |
| | 31 | -0.377000 | 0.640000 | -0.640000 | 0.000000 | 0.026532 | 0.000000 |

## Appendix II A Reference Code to Calculate $G(k;X)$

```
#include <stdio.h>
#include <stdlib.h>
#include <malloc.h>
#include <string.h>
#include <math.h>

int main()
{
	int test (unsigned int, unsigned int);

	printf("%d\n",test (5, 565));
	return 0;
}

int test (unsigned int kth, unsigned int len)
{
	printf("%d,%d\n", kth,len);

	int	s;
	unsigned __int64	i, j, k, r, m, h, b;
	unsigned __int64	n, x, t, d, t1;

	n = (unsigned __int64)pow(10, 7);
	n = len * n;

	constexpr	int	block = 512;
	constexpr	unsigned __int64 I = 0xffffffffffffffff;

	unsigned __int64 * pt;
	unsigned __int64 * pt2;
	pt = (unsigned __int64 *)calloc(n + 4, sizeof(unsigned __int64));
	pt2 = (unsigned __int64 *)calloc(n + 4, sizeof(unsigned __int64));

	if (pt == NULL)
	{
		printf("calloc is failed 1");
		exit(-1);
	}
	if (pt2 == NULL)
	{
		printf("calloc is failed 2");
```

```
        exit(-1);
    }

    pt[n] = 2, pt2[n] = 2;

    double w, v;
    v = 64.0;
    v = v * n;

    w = pow(v, 1.0 / kth);
    m = (unsigned __int64)ceil(w);

    unsigned __int64 *a;
    a = (unsigned __int64 *)calloc(m + 1, sizeof(unsigned __int64));
    if (a == NULL)
    {
        printf("calloc is failed 3");
        exit(-1);
    }

    unsigned __int64 id;
    id = 1;
    for (i = 0; i < m; i++)
    {
        x = i;
        for (r = 1; r < kth; r++)
        x = x * i;
        a[i] = x;
        t = x & 63;
        d = x >> 6;
        pt[d] ^= (id << t);
        pt2[d] ^= (id << t);
    }

    unsigned __int64 *bpt, *cpt;
    unsigned __int64 *jpt;

    h = (n >> 1);

    for (s = 2; s <= 100; s++)
    {
        printf("%d\n", s);
        j = 0;
        do
```

```
{
        cpt = pt + j;
        for ( k=j; *(cpt) == 0; k++ )
        {
                cpt++;
        }

        for (i = 1; i < m; i++)
        {

                x = a[i];
                t = x & 63;
                d = x >> 6;
                d = d + k;
                if (d > n)
                        break;
                else
                {
                        bpt = cpt;
                        jpt = pt2 + d;
                        if (d + block > n)
                                b = n - d + 1;
                        else
                                b = block;
                        if (t == 0)
                        {
                                for (r = 0; r < b; r++)
                                {
                                        *(jpt) |= *(bpt);
                                        jpt++, bpt++;
                                }
                        }
                        else
                        {
                                t1 = 64 - t;
                                x = *(bpt);
                                *(jpt) |= (x << t);
                                for (r = 1; r < b; r++)
                                {
                                        jpt++, bpt++;
                                        *(jpt) |= ((x >> t1) | (*(bpt)<< t));
                                        x = *(bpt);
                                }
                                *(jpt + 1) |= (x >> t1);
```

```
                    }
                }

            }
            j = k + block;
        } while (j <= n);

        memcpy(pt, pt2, sizeof(unsigned __int64)*(n+1));

        for (i = (n - 1); (pt2[i] == I) && (i >= h); i--)
        {
        }

        if (i < h)
        break;
    }

    free(pt2);
    free(pt);
    free(a);
    return(s);
}
```